\documentclass[reqno]{amsart}
\usepackage{amssymb}
\usepackage[vmarginratio=1:1,hmarginratio=1:1,body={16cm,23cm}]{geometry}%
\usepackage{ifpdf}
\ifpdf
 \usepackage[hyperindex]{hyperref}%,pagebackref
\else
 \expandafter\ifx\csname dvipdfm\endcsname\relax
 \usepackage[hypertex,hyperindex]{hyperref}
 \else
 \usepackage[dvipdfm,hyperindex]{hyperref}
 \fi
\fi
\theoremstyle{plain}
\newtheorem{theorem}{Theorem}[section]

\theoremstyle{remark}
\newtheorem{remark}{Remark}[section]
\numberwithin{equation}{section}
\DeclareMathOperator{\td}{d}
\DeclareMathOperator{\sech}{sech}
\DeclareMathOperator{\csch}{csch}
\newcommand{\bell}{\textup{B}}
\allowdisplaybreaks[4]

\begin{document}

\title[Derivatives of tangent and applications]
{Derivatives of tangent function and tangent numbers}

\author[F. Qi]{Feng Qi}
\address{Institute of Mathematics, Henan Polytechnic University, Jiaozuo City, Henan Province, 454010, China; College of Mathematics, Inner Mongolia University for Nationalities, Tongliao City, Inner Mongolia Autonomous Region, 028043, China; Department of Mathematics, College of Science, Tianjin Polytechnic University, Tianjin City, 300387, China}
\email{\href{mailto: F. Qi <qifeng618@gmail.com>}{qifeng618@gmail.com}, \href{mailto: F. Qi <qifeng618@hotmail.com>}{qifeng618@hotmail.com}, \href{mailto: F. Qi <qifeng618@qq.com>}{qifeng618@qq.com}}
\urladdr{\url{http://qifeng618.wordpress.com}}

\begin{abstract}
In the paper, by induction, the Fa\`a di Bruno formula, and some techniques in the theory of complex functions, the author finds explicit formulas for higher order derivatives of the tangent and cotangent functions as well as powers of the sine and cosine functions, obtains explicit formulas for two Bell polynomials of the second kind for successive derivatives of sine and cosine functions, presents curious identities for the sine function, discovers explicit formulas and recurrence relations for the tangent numbers, the Bernoulli numbers, the Genocchi numbers, special values of the Euler polynomials at zero, and special values of the Riemann zeta function at even numbers, and comments on five different forms of higher order derivatives for the tangent function and on derivative polynomials of the tangent, cotangent, secant, cosecant, hyperbolic tangent, and hyperbolic cotangent functions.
\end{abstract}

\keywords{explicit formula; recurrence; derivative; derivative polynomial; tangent; cotangent; Riemann zeta function; Bernoulli number; Genocchi number; Bell polynomial of the second kind; Euler polynomial; induction; Fa\`a di Bruno formula}

\subjclass[2010]{Primary 11B68, 11B83, 33B10; Secondary 11C08, 11M06, 26A24}

\thanks{This paper was typeset using \AmS-\LaTeX}

\maketitle

\section{Main results}

It is well known that the tangent function $\tan x$ can be expanded into the Maclaurin series
\begin{equation*}
\tan x=\sum_{k=1}^\infty\frac{(-1)^{k-1}2^{2k}\bigl(2^{2k}-1\bigr) B_{2k}}{(2k)!}x^{2k-1}
=\sum_{k=1}^\infty T_{2k-1}\frac{x^{2k-1}}{(2k-1)!}, \quad|x|<\frac{\pi}2,
\end{equation*}
see~\cite[p.~75, 4.3.67]{abram} and~\cite[p.~259]{Comtet-Combinatorics-74}, where $T_{2k-1}$ are called the tangent numbers or zag numbers and $B_{n}$ for $n\ge0$ are the Bernoulli numbers which may be defined by the power series expansion
\begin{equation*}
\frac{x}{e^x-1}=\sum_{n=0}^\infty B_n\frac{x^n}{n!}=1-\frac{x}2+\sum_{k=1}^\infty B_{2k}\frac{x^{2k}}{(2k)!}, \quad \vert x\vert<2\pi.
\end{equation*}
\par
The tangent numbers $T_{2k-1}$ may also be defined combinatorially as the numbers of alternating permutations on $2k-1=1,3,5,7,\dotsc$ symbols (where permutations that are the reverses of one another counted as equivalent). The first few tangent numbers $T_{2k-1}$ for $k=1,2,\dotsc,5$ are $1,2,16,272,7936$.
\par
It is clear that
\begin{equation}\label{tan-Bern-lim-der}
T_{2k-1}=(-1)^{k-1}\frac{2^{2k-1}\bigl(2^{2k}-1\bigr)}{k} B_{2k}=\lim_{x\to0}\tan^{(2k-1)}x.
\end{equation}
Consequently, one way to compute the tangent numbers $T_{2k-1}$ and the Bernoulli numbers $B_{2k}$ is to find explicit formulas for $\tan^{(2k-1)}x$.
To the best of the author's ability, explicit formulas for the $n$-th derivatives $(\tan x)^{(n)}$ of the tangent function $\tan x$ can not be searched out anywhere.
\par
It is well known that the Riemann zeta function $\zeta(s)$ may be defined by
\begin{equation*}
\zeta(s)=\sum_{k=1}^\infty\frac1{k^s},\quad \Re(s)>0,
\end{equation*}
that the Genocchi numbers $G_k$ are given by the generating function
\begin{equation*}
\frac{2z}{e^z+1}=\sum_{k=1}^\infty G_k\frac{z^k}{k!}, \quad |z|<\pi,
\end{equation*}
and that the Euler polynomials $E_k(x)$ are defined by
\begin{equation*}
\frac{2e^{xz}}{e^z+1}=\sum_{k=0}^\infty E_k(x)\frac{z^k}{k!}, \quad |z|<\pi.
\end{equation*}
The tangent numbers $T_{2k-1}$, the Bernoulli numbers $B_{2k}$, the Genocchi numbers $G_k$, the Euler polynomials $E_k(x)$, and the Riemann zeta function $\zeta(s)$ have closely relations.
\par
In combinatorics, the Bell polynomials of the second kind, or say, the partial Bell polynomials, denoted by $\bell_{n,k}(x_1,x_2,\dotsc,x_{n-k+1})$ for $n\ge k\ge0$, are defined by
\begin{equation*}
\bell_{n,k}(x_1,x_2,\dotsc,x_{n-k+1})=\sum_{\substack{1\le q\le n,\ell_q\in\{0\}\cup\mathbb{N}\\ \sum_{q=1}^ni\ell_q=n\\
\sum_{q=1}^n\ell_q=k}}\frac{n!}{\prod_{q=1}^{n-k+1}\ell_q!} \prod_{q=1}^{n-k+1}\Bigl(\frac{x_q}{q!}\Bigr)^{\ell_q}.
\end{equation*}
See~\cite[p.~134, Theorem~A]{Comtet-Combinatorics-74}. In combinatorial analysis, the Fa\`a di Bruno formula plays an important role and may be described in terms of the Bell polynomials of the second kind $\bell_{n,k}$ by
\begin{equation}\label{Bruno-Bell-Polynomial}
\frac{\td^n}{\td t^n}f\circ h(t)=\sum_{k=1}^nf^{(k)}(h(t)) \bell_{n,k}\bigl(h'(t),h''(t),\dotsc,h^{(n-k+1)}(t)\bigr).
\end{equation}
See~\cite[p.~139, Theorem~C]{Comtet-Combinatorics-74}.
\par
In this paper, by induction, the Fa\`a di Bruno formula, and some techniques in the theory of complex functions, we will establish general and explicit formulas for the $n$-th derivatives of $\tan x$, $\cot x$, $\sin^kx$, and $\cos^kx$ for $k\in\mathbb{N}$, presents curious identities for the sine function, and obtain explicit formulas for two Bell polynomials of the second kind
\begin{equation*}
\bell_{m,k}\biggl(-\sin x,-\cos x,\sin x,\cos x,\dotsc, -\sin\biggl[x+(m-k)\frac{\pi}{2}\biggr]\biggr)
\end{equation*}
and
\begin{equation*}
\bell_{m,k}\biggl(\cos x,-\sin x,-\cos x,\sin x,\dotsc, -\cos\biggl[x+(m-k)\frac{\pi}{2}\biggr]\biggr).
\end{equation*}
By applying these formulas, we will also derive explicit formulas and recurrence relations for the tangent numbers $T_{2n-1}$, the Bernoulli numbers $B_{2n}$, the Genocchi numbers $G_{2n}$, special values $E_{2n-1}(0)$ of the Euler polynomials at $0$, and special values $\zeta(2n)$ of the Riemann zeta function $\zeta(z)$ at even numbers $2n$. Finally, we will comment on five different forms of higher order derivatives for $\tan x$ and on derivative polynomials of $\tan x$, $\cot x$, $\sec x$, $\csc x$, $\tanh x$, and $\coth x$.
\par
Our main results may be stated as the following theorems.

\begin{theorem}\label{tan-deriv-thm}
For $n\in\mathbb{N}$, derivatives of the tangent and cotangent functions may be computed by
\begin{multline}\label{tan-deriv-unity}
\tan^{(n)}x=\frac1{\cos^{n+1}x}\Biggl\{\frac12\biggl[1+\frac{1+(-1)^n}2\biggr]a_{n,\frac{1+(-1)^n}2} \sin\biggl[\frac{1+(-1)^n}2x +\frac{1-(-1)^n}2\frac\pi2\biggr]\\
+\sum_{k=1}^{\frac12\bigl[n-1-\frac{1+(-1)^n}2\bigr]} a_{n,2k+\frac{1+(-1)^n}2}\sin\biggl[\biggl(2k+\frac{1+(-1)^n}2\biggr)x +\frac{1-(-1)^n}2\frac\pi2\biggr]\Biggr\}
\end{multline}
and
\begin{multline}\label{cotan-deriv-n}
\cot^{(n)}x=\frac{(-1)^n}{\sin^{n+1}x}\Biggl\{\frac{1}{2} \biggl[1+\frac{1+(-1)^n}{2}\biggr]a_{n,\frac{1+(-1)^n}2} \cos\biggl[\frac{1+(-1)^n}2x\biggr]\\*
+\sum_{k=1}^{\frac12\bigl[n-1-\frac{1+(-1)^n}2\bigr]} (-1)^ka_{n,2k+\frac{1+(-1)^n}2} \cos\biggl[\biggl(2k+\frac{1+(-1)^n}2\biggr)x\biggr]\Biggr\},
\end{multline}
where
\begin{equation}\label{a(p-q)-unify}
a_{p,q}=(-1)^{\frac12\bigl[q-\frac{1+(-1)^p}2\bigr]}[1-(-1)^{p-q}]\sum_{\ell=0}^{\frac{p-q-1}2} (-1)^{\ell}\binom{p+1}{\ell}\biggl(\frac{p-q-1}2-\ell+1\biggr)^p
\end{equation}
for $p>q\ge0$.
\end{theorem}

\begin{theorem}\label{bell-sin-cos-thm}
Let $k,m$ be nonnegative integers such that $(k,m)\ne(0,0)$.
\par
For all $k,m\ge0$, the sine and cosine functions satisfy
\begin{gather}\label{sin-ell-d}
\frac{\td^m\sin^kx}{\td x^m}
=\frac{(-1)^k}{2^k} \sum_{q=0}^k(-1)^q\binom{k}{q}(2q-k)^m \cos\biggl[(m-k)\frac\pi2+(2q-k)x\biggr],\\
\label{cos-ell-derr}
\frac{\td^m\cos^kx}{\td x^m}
=\frac1{2^k}\sum_{q=0}^k\binom{k}{q}(2q-k)^m \cos\biggl[\frac{\pi}2m+(2q-k)x\biggr],\\
\label{eq=0-2-sin-cos}
\sum_{q=0}^k(-1)^q\binom{k}{q}(2q-k)^m \sin\biggl[(m-k)\frac\pi2+(2q-k)x\biggr]=0,
\end{gather}
and
\begin{equation}\label{eq=0-1-sin-cos}
\sum_{q=0}^k\binom{k}{q}(2q-k)^m \sin\biggl[\frac{\pi}2m+(2q-k)x\biggr]=0.
\end{equation}
\par
For $m\ge k\ge1$, the Bell polynomials of the second kind $\bell_{m,k}$ satisfy
\begin{multline}\label{bell-sin-eq}
\bell_{m,k}\biggl(-\sin x,-\cos x,\sin x,\cos x,\dotsc, -\sin\biggl[x+(m-k)\frac{\pi}{2}\biggr]\biggr)\\
=\frac{(-1)^k}{k!}\cos^kx\sum_{\ell=0}^k\frac{(-1)^\ell}{2^\ell}\binom{k}{\ell}
\frac1{\cos^{\ell}x}\sum_{q=0}^\ell\binom{\ell}{q}(2q-\ell)^m \cos\biggl[\frac{\pi}2m+(2q-\ell)x\biggr],
\end{multline}
and
\begin{multline}\label{bell-sin=ans}
\bell_{m,k}\biggl(\cos x,-\sin x,-\cos x,\sin x,\dotsc, -\cos\biggl[x+(m-k)\frac{\pi}{2}\biggr]\biggr)\\
=\frac{(-1)^k}{k!}\sin^{k}x\sum_{\ell=0}^k\frac1{2^\ell}\binom{k}{\ell}\frac1{\sin^{\ell}x}
\sum_{q=0}^\ell(-1)^q\binom{\ell}{q}(2q-\ell)^m \cos\biggl[(m-\ell)\frac\pi2+(2q-\ell)x\biggr].
\end{multline}
\par
For $n\in\mathbb{N}$, derivatives of the tangent and cotangent functions may be computed by
\begin{equation}\label{tan-der-nth-eq}
\frac{\td^n\tan x}{\td x^n}=
-\sum_{k=1}^{n+1}\frac1{k}\sum_{\ell=0}^k\frac{(-1)^\ell}{2^\ell}\binom{k}{\ell}
\frac1{\cos^{\ell}x}\sum_{q=0}^\ell\binom{\ell}{q}(2q-\ell)^{n+1} \sin\biggl[\frac{\pi}2n+(2q-\ell)x\biggr]
\end{equation}
and
\begin{equation}\label{cotan-der-nth-eq}
\frac{\td^n\cot x}{\td x^n}
=\sum_{k=1}^{n+1}\frac1k\sum_{\ell=0}^k\frac1{2^\ell}\binom{k}{\ell}\frac1{\sin^{\ell}x}
\sum_{q=0}^\ell(-1)^q\binom{\ell}{q}(2q-\ell)^{n+1} \sin\biggl[(n-\ell)\frac\pi2+(2q-\ell)x\biggr].
\end{equation}
\end{theorem}

\begin{theorem}\label{tan-cot-app-Bernoulli}
For $m\in\mathbb{N}$, the tangent numbers $T_{2m-1}$ may be computed by
\begin{equation}\label{Tangent-No-eq}
T_{2m-1}=\sum_{k=0}^{m-1} (-1)^k\binom{2m}{k}(m-k)^{2m-1}
+2\sum_{k=1}^{m-1}(-1)^k\sum_{\ell=0}^{m-k-1} (-1)^{\ell}\binom{2m}{\ell}(m-k-\ell)^{2m-1}
\end{equation}
and
\begin{equation}\label{tan-No-faaa}
T_{2m-1}=(-1)^m\sum_{\ell=1}^{2m}(-1)^\ell\frac1{2^\ell} \biggl(\frac1\ell-\frac1{m+1}\biggr) \binom{2m+1}{\ell} \sum_{q=0}^\ell\binom{\ell}{q}(2q-\ell)^{2m}.
\end{equation}
and satisfy a recurrence
\begin{equation}\label{tan-No-rec-eq}
T_{2m+1}=\sum_{k=1}^m\binom{2m}{2k-1}T_{2k-1}T_{2(m-k)+1}.
\end{equation}
\end{theorem}

\begin{theorem}\label{Ber-zeta-genoc-thm}
For $m\in\mathbb{N}$, the Bernoulli numbers $B_{2m}$ may be computed by
\begin{multline}\label{Tangent-No-eq-B}
B_{2m}=(-1)^{m-1}\frac{m} {2^{2m-1}\bigl(2^{2m}-1\bigr)}\Biggl[\sum_{k=0}^{m-1} (-1)^k\binom{2m}{k}(m-k)^{2m-1}\\*
+2\sum_{k=1}^{m-1}(-1)^k\sum_{\ell=0}^{m-k-1} (-1)^{\ell}\binom{2m}{\ell}(m-k-\ell)^{2m-1}\Biggr]
\end{multline}
and
\begin{equation}\label{tan-No-faaa-B}
B_{2m}=-\frac{m} {2^{2m-1}\bigl(2^{2m}-1\bigr)}\sum_{\ell=1}^{2m}(-1)^\ell\frac1{2^\ell} \biggl(\frac1\ell-\frac1{m+1}\biggr) \binom{2m+1}{\ell} \sum_{q=0}^\ell\binom{\ell}{q}(2q-\ell)^{2m}
\end{equation}
and the sequence $\mathcal{B}_{2m}=\frac{2^{2m}-1}{2m} B_{2m}$ satisfies a recurrence
\begin{equation}\label{tan-No-rec-eq-B}
\mathcal{B}_{2(m+1)}=-\sum_{k=1}^m\binom{2m}{2k-1}\mathcal{B}_{2k}\mathcal{B}_{2(m-k+1)}.
\end{equation}
\par
For $m\in\mathbb{N}$, special values $\zeta(2m)$ of the zeta function $\zeta(z)$ at $z=2m$ may be computed by
\begin{multline}\label{zeta-No-eq-a}
\zeta(2m)=\frac{\pi^{2m}m}{(2m)!(2^{2m}-1)}\Biggl[\sum_{k=0}^{m-1} (-1)^k\binom{2m}{k}(m-k)^{2m-1}\\
+2\sum_{k=1}^{m-1}(-1)^k\sum_{\ell=0}^{m-k-1} (-1)^{\ell}\binom{2m}{\ell}(m-k-\ell)^{2m-1}\Biggr]
\end{multline}
and
\begin{equation}\label{zeta-No-eq-b}
\zeta(2m)=(-1)^m\frac{\pi^{2m}m}{(2m)!(2^{2m}-1)} \sum_{\ell=1}^{2m}(-1)^\ell\frac1{2^\ell} \biggl(\frac1\ell-\frac1{m+1}\biggr) \binom{2m+1}{\ell} \sum_{q=0}^\ell\binom{\ell}{q}(2q-\ell)^{2m}
\end{equation}
and the sequence $\mathcal{Z}_{2m}=(-1)^m\frac{(2m)!(2^{2m}-1)}{(2\pi)^{2m}m}\zeta(2m)$ satisfies a recurrence
\begin{equation}\label{zeta-No-rec-eq-B}
\mathcal{Z}_{2(m+1)}=\sum_{k=1}^m\binom{2m}{2k-1}\mathcal{Z}_{2k}\mathcal{Z}_{2(m-k+1)}.
\end{equation}
\par
For $m\in\mathbb{N}$, the Genocchi numbers $G_{2m}$ may be calculated by
\begin{multline}\label{G-No-eq-a}
G_{2m}=(-1)^m\frac{m} {2^{2(m-1)}}\Biggl[\sum_{k=0}^{m-1} (-1)^k\binom{2m}{k}(m-k)^{2m-1}\\
+2\sum_{k=1}^{m-1}(-1)^k\sum_{\ell=0}^{m-k-1} (-1)^{\ell}\binom{2m}{\ell}(m-k-\ell)^{2m-1}\Biggr]
\end{multline}
and
\begin{equation}\label{G-No-eq-b}
G_{2m}=\frac{2m} {2^{2m-1}}\sum_{\ell=1}^{2m}(-1)^\ell\frac1{2^\ell} \biggl(\frac1\ell-\frac1{m+1}\biggr) \binom{2m+1}{\ell} \sum_{q=0}^\ell\binom{\ell}{q}(2q-\ell)^{2m}
\end{equation}
and the sequence $\mathcal{G}_{2m}=\frac1{m}G_{2m}$ satisfies a recurrence
\begin{equation}\label{No-rec-eq-G}
\mathcal{G}_{2(m+1)}=\frac14\sum_{k=1}^m\binom{2m}{2k-1}\mathcal{G}_{2k}\mathcal{G}_{2(m-k+1)}.
\end{equation}
\par
For $m\in\mathbb{N}$, special values $E_{2m-1}(0)$ of the Euler polynomials $E_m(x)$ at $x=0$ may be calculated by
\begin{multline}\label{E-No-eq-a}
E_{2m-1}(0)=(-1)^m\frac1{2^{2m-1}}\Biggl[\sum_{k=0}^{m-1} (-1)^k\binom{2m}{k}(m-k)^{2m-1}\\
+2\sum_{k=1}^{m-1}(-1)^k\sum_{\ell=0}^{m-k-1} (-1)^{\ell}\binom{2m}{\ell}(m-k-\ell)^{2m-1}\Biggr]
\end{multline}
and
\begin{equation}\label{E-No-eq-b}
E_{2m-1}(0)=\frac1{2^{2m-1}} \sum_{\ell=1}^{2m}(-1)^\ell\frac1{2^\ell} \biggl(\frac1\ell-\frac1{m+1}\biggr) \binom{2m+1}{\ell} \sum_{q=0}^\ell\binom{\ell}{q}(2q-\ell)^{2m}
\end{equation}
and the sequence $E_{2m-1}(0)$ satisfies a recurrence
\begin{equation}\label{tan-No-rec-eq-E}
E_{2m+1}(0)=\frac12\sum_{k=1}^m\binom{2m}{2k-1}E_{2k-1}(0)E_{2(m-k)+1}(0).
\end{equation}
\end{theorem}

\section{Comparisons and highlights}

Before proving our main results, we compare them with some known conclusions for showing their highlights.

\subsection{Applications of the Bell polynomials~\eqref{bell-sin-eq} and~\eqref{bell-sin=ans}}

By the Fa\'a di Bruno formula~\eqref{Bruno-Bell-Polynomial} and the Bell polynomials~\eqref{bell-sin-eq} and~\eqref{bell-sin=ans} in Theorem~\ref{bell-sin-cos-thm}, it is easy to write
\begin{equation*}
(\sec x)^{(n)}=\biggl(\frac1{\cos x}\biggr)^{(n)}
=\frac{1}{\cos x}\sum_{k=1}^n \sum_{\ell=0}^k\frac{(-1)^\ell}{2^\ell}\binom{k}{\ell}
\frac1{\cos^{\ell}x}\sum_{q=0}^\ell\binom{\ell}{q}(2q-\ell)^n \cos\biggl[\frac{\pi}2n+(2q-\ell)x\biggr]
\end{equation*}
and
\begin{equation*}
(\csc x)^{(n)}=\biggl(\frac1{\sin x}\biggr)^{(n)}
=\frac1{\sin x}\sum_{k=1}^n\sum_{\ell=0}^k\frac1{2^\ell}\binom{k}{\ell}\frac1{\sin^{\ell}x}
\sum_{q=0}^\ell(-1)^q\binom{\ell}{q}(2q-\ell)^n \cos\biggl[(n-\ell)\frac\pi2+(2q-\ell)x\biggr].
\end{equation*}
Consequently, we may find explicit formulas for the secant and cosecant numbers. Due to the limitation of length, we do not write down them in details.
\par
The formulas~\eqref{bell-sin-eq} and~\eqref{bell-sin=ans} answer a problem in~\cite[Section~5]{Deriv-Arcs-Cos.tex}. These formulas, together with the Fa\`a di Bruno formula~\eqref{Bruno-Bell-Polynomial}, may be applied to calculate the $n$-th derivatives of functions of the forms $f(\sin x)$ and $f(\cos x)$, such as $(\sin x)^\alpha$, $(\cos x)^\alpha$, $(\sec x)^\alpha$, $(\csc x)^\alpha$, $e^{\pm\sin x}$, $e^{\pm\cos x}$, $\ln(\cos x)$, $\ln(\sin x)$, $\ln(\sec x)$, $\ln(\csc x)$, $\sin(\sin x)$, $\cos(\sin x)$, $\sin(\cos x)$, and $\cos(\cos x)$.

\subsection{Five different forms of higher order derivatives for $\boldsymbol{\tan x}$}

On~\cite[pp.~28--31, Chapter~II]{Schwatt-1924},
see also its second edition~\cite{Schwatt-1962}, by virtue of the formula
\begin{equation*}
\frac{\td^ny}{\td x^n}=\sum_{k=1}^n\frac{(-1)^k}{k!} \sum_{\alpha=1}^k(-1)^\alpha\binom{k}{\alpha}u^{k-\alpha}
\frac{\td^n(u^\alpha)}{\td x^n}\frac{\td^ky}{\td u^k}
\end{equation*}
in~\cite[p.~12, (83)]{Schwatt-1924}, where $y=\phi(u)$ and $u=f(x)$, the formulas
\begin{align}\label{Schwatt-1924-p29-1}
\frac{\td^{2n}\tan x}{\td x^{2n}}&=(-1)^{n+1}2^{2n}\sum_{k=1}^{2n}\frac1{2^k}\sec^{k+1}x\sin[(k-1)x] \sum_{\alpha=1}^k(-1)^\alpha\binom{k}{\alpha}\alpha^{2n},\\
\frac{\td^{2n+1}\tan x}{\td x^{2n+1}} &=(-1)^{n+1}2^{2n+1}\sum_{k=1}^{2n+1}\frac1{2^k}\sec^{k+1}x\cos[(k-1)x] \sum_{\alpha=1}^k(-1)^\alpha\binom{k}{\alpha}\alpha^{2n+1},\label{Schwatt-1924-p29-2}\\
\label{Schwatt-1924-p29-comb}
\frac{\td^n\tan x}{\td x^n} &=(-1)^{\frac12\bigl[n+\frac{1-(-1)^n}2\bigr]}2^n\sum_{k=1}^n\frac1{2^k}\sec^{k+1}x \sin\biggl[\frac\pi2\frac{1-(-1)^n}2+(k-1)x\biggr] \sum_{\alpha=1}^k(-1)^\alpha\binom{k}{\alpha}\alpha^n\\
\frac{\td^{2n}\tan x}{\td x^{2n}}&=(-1)^{n-1}2^{2n}\sec^2x\sum_{k=1}^{2n}\frac1{2^k} \sum_{\alpha=1}^k(-1)^\alpha\binom{k}{\alpha}\alpha^{2n}N_{2\beta+1},\\
\frac{\td^{2n+1}\tan x}{\td x^{2n+1}}&=(-1)^{n-1}2^{2n+1}\sec^2x\sum_{k=1}^{2n+1}\frac1{2^k} \sum_{\alpha=1}^k(-1)^\alpha\binom{k}{\alpha}\alpha^{2n+1}N_{2\beta},\\
\frac{\td^{n}\tan x}{\td x^{n}}&=(-1)^{\lfloor\frac{n+2}2\rfloor} 2^{n}\sec^2x\sum_{k=1}^{n}\frac1{2^k} \sum_{\alpha=1}^k(-1)^\alpha\binom{k}{\alpha}\alpha^{n}N_{2\beta+\frac{1+(-1)^n}2} \label{Schwatt-1924-comb}
\end{align}
were obtained for $n\in\mathbb{N}$, where
\begin{equation*}
N_{2\beta}=\sum_{\beta=0}^{\lfloor\frac{k-1}2\rfloor}(-1)^\beta\binom{k-1}{2\beta}\tan^{2\beta}x,
\quad
N_{2\beta+1}=\sum_{\beta=0}^{\lfloor\frac{k-2}2\rfloor}(-1)^\beta\binom{k-1}{2\beta+1}\tan^{2\beta+1}x,
\end{equation*}
and $\lfloor x\rfloor$, whose value is the biggest integer not more than $x$, stands for the floor function of $x$.
\par
The formula~\eqref{tan-deriv-unity} may be separately written as
\begin{equation}\label{tan-deriv-2n-1}
\tan^{(2n-1)}x=\frac1{\cos^{2n}x}\sum_{k=0}^{n-1}a_{2n-1,2k}\cos(2kx)
\end{equation}
and
\begin{equation}\label{tan-deriv-2n}
\tan^{(2n)}x=\frac1{\cos^{2n+1}x}\sum_{k=0}^{n-1}a_{2n,2k+1}\sin[(2k+1)x]
\end{equation}
for $n\in\mathbb{N}$, where
\begin{align*}
a_{1,0}&=1,\\
a_{2n-1,0}&=2n\sum_{\ell=0}^{n-2}(-1)^{\ell}\binom{2n-1}{\ell}(n-\ell-1)^{2n-2}
\end{align*}
for $n>1$, and
\begin{equation*}
a_{p,q}=(-1)^{\frac12\bigl[p-\frac{3+(-1)^p}2\bigr]}2\sum_{\ell=0}^{\frac{p-q-1}2} (-1)^{\frac{p-q-1}2-\ell}\binom{p+1}{\ell}\biggl(\frac{p-q-1}2-\ell+1\biggr)^p
\end{equation*}
for $0<q<p$ with $p-q$ being a positive odd number. See the first version of the preprint~\cite{derivative-tan-cot.tex}. It is clear that, to some extent, the forms of the formulas~\eqref{tan-deriv-2n-1} and~\eqref{tan-deriv-2n} are apparently simpler than those of the formulas~\eqref{Schwatt-1924-p29-1} and~\eqref{Schwatt-1924-p29-2}. On the other hand, there are $2n$ and $2n+1$ terms in the formulas~\eqref{Schwatt-1924-p29-1} and~\eqref{Schwatt-1924-p29-2} respectively, but only $n$ terms in both of the formulas~\eqref{tan-deriv-2n-1} and~\eqref{tan-deriv-2n}. This means that the formula~\eqref{tan-deriv-unity} in Theorem~\ref{tan-deriv-thm} is essentially simpler than~\eqref{Schwatt-1924-p29-comb}.
\par
Let $u=u(x)$ and $v=v(x)\ne0$ be differentiable functions. In~\cite[p.~40]{Bourbaki-Spain-2004}, the formula
\begin{equation}\label{Sitnik-Bourbaki}
\frac{\td^n}{\td x^n}\biggl(\frac{u}{v}\biggr)
=\frac{(-1)^n}{v^{n+1}}
\begin{vmatrix}
u & v & 0 & \dots & 0\\
u' & v' & v & \dots & 0\\
u'' & v'' & 2v' & \dots & 0\\
\hdotsfor[2]{5}\\
u^{(n-1)} & v^{(n-1)} & \binom{n-1}1v^{(n-2)} &  \dots & v\\
u^{(n)} & v^{(n)} & \binom{n}1v^{(n-1)} & \dots & \binom{n}{n-1}v'
\end{vmatrix}
\end{equation}
for the $n$th derivative of the ratio $\frac{u(x)}{v(x)}$ was listed. For easy understanding and convenient availability, we now reformulate the formula~\eqref{Sitnik-Bourbaki} as
\begin{equation}\label{Sitnik-Bourbaki-reform}
\frac{\td^n}{\td x^n}\biggl(\frac{u}{v}\biggr)
=\frac{(-1)^n}{v^{n+1}}
\begin{vmatrix}
A_{(n+1)\times1}&B_{(n+1)\times n}
\end{vmatrix}_{(n+1)\times(n+1)},
\end{equation}
where $|\cdot|_{(n+1)\times(n+1)}$ denotes a determinant and the matrices
\begin{equation*}
A_{(n+1)\times1}=(a_{\ell,1})_{0\le \ell\le n}
\end{equation*}
and
\begin{equation*}
B_{(n+1)\times n}=(b_{\ell,m})_{0\le \ell\le n,0\le m\le n-1}
\end{equation*}
satisfy
\begin{equation*}
a_{\ell,1}=u^{(\ell)}(x)\quad \text{and}\quad b_{\ell,m}=\binom{\ell}{m}v^{(\ell-m)}(x)
\end{equation*}
under the conventions that $v^{(0)}(x)=v(x)$ and that $\binom{p}{q}=0$ and $v^{(p-q)}(x)\equiv0$ for $p<q$. See also~\cite[Lemma~2.1]{Euler-No-3Sum.tex}. Applying the formula~\eqref{Sitnik-Bourbaki-reform} to $\tan x=\frac{\sin x}{\cos x}$ acquires
\begin{equation}\label{a-b-relation-ratio}
a_{\ell,1}=\sin^{(\ell)}x=\sin\biggl(x+\frac\pi2\ell\biggr)\quad\text{and}\quad
b_{\ell,m}=\binom{\ell}{m}\cos\biggl(x+\frac\pi2(\ell-m)\biggr).
\end{equation}
Therefore, we can find an alternative form, the fifth form, for higher order derivatives of $\tan x$. The first four forms for higher order derivatives of $\tan x$ are~\eqref{tan-deriv-unity}, \eqref{tan-der-nth-eq}, \eqref{Schwatt-1924-p29-comb}, and~\eqref{Schwatt-1924-comb}. These five different forms come from the induction and different formulas for higher order derivatives of composite functions.

\subsection{Derivative polynomials}

Suppose $f$ is a function whose derivative is a polynomial in $f$, that is, $f'(x)=P(f(x))$ for some polynomial $P$. Then all the higher order derivatives of $f$ are also polynomials in $f$, so we have a sequence of polynomials $P_n$ defined by $f^{(n)}(x)=P_n(f(x))$ for $n\ge0$. As usual, we call $P_n(u)$ the derivative polynomials of $f$.
\par
In terms of this terminology, some results in~\cite{exp-derivative-sum-Combined.tex, Eight-Identy-More.tex-JCAM, CAM-D-13-01430-Xu-Cen} may be restated as follows: when $\lambda>0$ and $t\ne-\frac{\ln\lambda}\alpha$ or when $\lambda<0$ and $t\in\mathbb{R}$, the derivative polynomials of the function $\frac1{\lambda e^{\alpha t}-1}$ are
\begin{equation}\label{der-polynomial-exp}
(-1)^n\alpha^n\sum_{m=1}^{n+1}{(m-1)!S(n+1,m)}u^m
\end{equation}
in the variable $u$ for $\alpha\ne0$ and $n\in\mathbb{N}$, where
\begin{equation*}
S(q,m)=\frac1{m!}\sum_{\ell=1}^m(-1)^{m-\ell}\binom{m}{\ell}\ell^{q}
\end{equation*}
for $1\le m\le q$ are the Stirling numbers of the second kind which may be generated by
\begin{equation*}
\frac{(e^x-1)^k}{k!}=\sum_{n=k}^\infty S(n,k)\frac{x^n}{n!}, \quad k\in\mathbb{N}.
\end{equation*}
\par
In~\cite[p.~25, (5)]{Hoffman-Monthly1995} and~\cite[Theorem~3.2]{Hoffman-Monthly1995}, it was obtained that the derivative polynomials $P_n$ and $Q_n$ defined by
\begin{equation*}
\frac{\td^n(\tan x)}{\td x^n}=P_n(\tan x)\quad\text{and}\quad
\frac{\td^n(\sec x)}{\td x^n}=Q_n(\tan x)\sec x
\end{equation*}
for $n\ge0$ are polynomials of degree $n+1$ and $n$ respectively and satisfy the recurrences
\begin{equation}\label{recur-der-polyn}
\begin{aligned}
P_{n+1}(u)&=\sum_{k=0}^n\binom{n}{k}P_k(u)P_{n-k}(u)+\delta_{0n},\\
Q_{n+1}(u)&=\sum_{k=0}^n\binom{n}{k}P_k(u)Q_{n-k}(u),\\
P_{n+1}(u)&=\bigl(1+u^2\bigr)\sum_{k=0}^n\binom{n}{k}Q_k(u)Q_{n-k}(u),
\end{aligned}
\end{equation}
where
\begin{equation*}
\delta_{ij}=\begin{cases}
0,&i\ne j,\\ 1, & i=j,
\end{cases}
\quad P_0(u)=u, \quad\text{and}\quad P_1(u)=1+u^2.
\end{equation*}
In~\cite{Hoffman-EJC1999}, more recurrences for the derivative polynomials $P_n$ and $Q_n$ were obtained and they were applied to combinatorics and number theory.
We observe that the formulas
\begin{gather}\label{tan-deriv-Quadratic-form}
\tan^{(n+1)}x=\bigl(\tan^2x\bigr)^{(n)}=\sum_{k=0}^{n}\binom{n}{k}\tan^{(k)}x\tan^{(n-k)}x,\\
\begin{gathered}\notag
\sec^{(n+1)}x=(\sec x\tan x)^{(n)}=\sum_{k=0}^n\binom{n}{k}\sec^{(k)}x\tan^{(n-k)}x,\\
\tan^{(n+1)}x=\bigl(\sec^2x\bigr)^{(n)}=\sum_{k=0}^{n}\binom{n}{k}\sec^{(k)}x\sec^{(n-k)}x
\end{gathered}
\end{gather}
may be used to straightforwardly and simply recover the formulas in~\eqref{recur-der-polyn} for $n\in\mathbb{N}$.
\par
Let
\begin{equation*}
D=\frac{\td}{\td x}, \quad y=\tan x,\quad \text{and}\quad z=\sec x.
\end{equation*}
For $n\ge0$, define
\begin{equation}
\begin{gathered}\label{DY-power-1}
(Dy)^0(y)=y, \quad (Dy)(y)=D(y^2), \quad (Dy)^0(z)=z, \quad (Dy)(z)=D(yz),\\
(Dy)^{n+1}(y)=(Dy)(Dy)^n(y)=D\bigl(y(Dy)^n(y)\bigr),\\
(Dy)^{n+1}(z)=(Dy)(Dy)^n(z)=D\bigl(y(Dy)^n(z)\bigr),\\
(yD)^{n+1}(y)=(yD)(yD)^n(y)=yD\bigl((yD)^n(y)\bigr),\\
(yD)^{n+1}(z)=(yD)(yD)^n(z)=yD\bigl((yD)^n(z)\bigr).
\end{gathered}
\end{equation}
In~\cite{Ma-EJC-2013, Ma-BAMS-2014}, the above quantities were computed in terms of polynomials in variables $y$ and $z$ and connected with the Eulerian numbers and polynomials and others.
We observe the following contradiction: since
\begin{equation*}
(Dy)^2(y)=(Dy)(Dy)(y)=(Dy)D\big(y^2\bigr)
\end{equation*}
and
\begin{equation*}
(Dy)^2(y)=D\bigl(y(Dy)(y)\bigr)=D\bigl(yD\bigl(y^2\bigr)\bigr)=(Dy)D\bigl(y^2\bigr)+yD^2\bigl(y^2\bigr),
\end{equation*}
it follows that
\begin{equation*}
(Dy)D\big(y^2\bigr)=(Dy)D\bigl(y^2\bigr)+yD^2\bigl(y^2\bigr),
\end{equation*}
that is,
\begin{equation*}
yD^2\bigl(y^2\bigr)=0,
\end{equation*}
which means that $y^2=ax+b$, where $a,b\in\mathbb{C}$. Therefore, the authors of the papers~\cite{Ma-EJC-2013, Ma-BAMS-2014, Ma-Mansour-Wang-1403} and related ones using the definitions in~\eqref{DY-power-1} and their analogues should select a better manner and choose suitable symbols to express their ideas. The first author of the papers~\cite{Ma-EJC-2013, Ma-BAMS-2014, Ma-Mansour-Wang-1403} told the current author on May 24, 2015 that he explained in~\cite[Theorem~10]{MA-EJC} the above easily-confused definitions more explicitly.
\par
In~\cite{Boyadzhiev-Fibonacci2007, Boyadzhiev-arXiv0903}, the formulas
\begin{gather*}
\biggl(\frac{\td}{\td x}\biggr)^m\coth x=C_m(\coth x),\quad
\biggl(\frac{\td}{\td x}\biggr)^m\tanh x=C_m(\tanh x),\\
\biggl(\frac{\td}{\td x}\biggr)^m\sech x=(\sech x) S_m(\tanh x),\quad
\biggl(\frac{\td}{\td x}\biggr)^m\csch x=(\csch x) S_m(\coth x),\\
P_m(z)=i^{m+1}2^m(1-iz)\omega_m\biggl(-\frac{1+iz}{2}\biggr)
=-i^{m+1}(-2)^m(iz+1)\sum_{k=0}^m\frac{k!}{2^k}S(m,k)(iz-1)^k,\\
Q_m(z)=i^mS_m(iz)
\end{gather*}
were established for $m\ge1$, where $C_0(z)=1$,
\begin{gather*}
C_m(z)=(-2)^m(z+1)\omega_m\biggl(\frac{z-1}2\biggr)
=(-2)^m(z+1)\sum_{k=0}^m\frac{k!}{2^k}S(m,k)(z-1)^k,\quad m\ge1,\\
S_m(z)=\sum_{k=0}^m\binom{m}{k}2^k\omega_k\biggl(-\frac{z+1}2\biggr)
=\sum_{j=0}^m\Biggl[(-1)^jj!\sum_{k=j}^m\binom{m}kS(k,j)2^{k-j}\Biggr](z+1)^j,\quad m\ge0,
\end{gather*}
and
\begin{equation*}
\omega_n(x)=\sum_{k=0}^nS(n,k)k!x^k
\end{equation*}
are called the geometric polynomials. We notice that the formulas
\begin{gather*}
\tanh^{(n+1)}x=-\bigl(\tanh^2x\bigr)^{(n)}
=-\sum_{k=0}^n\binom{n}{k}\tanh^{(k)}x \tanh^{(n-k)}x,\\
\tanh^{(n+1)}x=\bigl(\sech^2x\bigr)^{(n)}
=\sum_{k=0}^n\binom{n}{k}\sech^{(k)}x \sech^{(n-k)}x,
\end{gather*}
and
\begin{equation*}
\sech^{(n+1)}x=-(\sech x\tanh x)^{(n)}=-\sum_{k=0}^n\sech^{(k)}x\tanh^{(n-k)}x
\end{equation*}
for $n\ge1$ imply trivially the recurrences
\begin{equation}
\begin{aligned}\label{recur-der-pol2}
C_{n+1}(u)&=-\sum_{k=0}^n\binom{n}{k}C_k(u)C_{n-k}(u),\\
C_{n+1}(u)&=\bigl(1+u^2\bigr)\sum_{k=0}^n\binom{n}{k}S_k(u)S_{n-k}(u),\\
S_{n+1}(u)&=-\sum_{k=0}^n\binom{n}{k}C_k(u)S_{n-k}(u)
\end{aligned}
\end{equation}
for $n\in\mathbb{N}$.
\par
Because
\begin{gather*}
\cot^{(n+1)}x =-\bigl(\cot^2x\bigr)^{(n)}=-\sum_{k=0}^{n}\binom{n}{k}\cot^{(k)}x \cot^{(n-k)}x,\\
\cot^{(n+1)}x =-\bigl(\csc^2x\bigr)^{(n)}=-\sum_{k=0}^{n}\binom{n}{k}\csc^{(k)}x \csc^{(n-k)}x,
\end{gather*}
and
\begin{equation*}
\csc^{(n+1)}x =-(\cot x \csc x)^{(n)}=-\sum_{k=0}^{n}\binom{n}{k}\cot^{(k)}x \csc^{(n-k)}x
\end{equation*}
for $n\in\mathbb{N}$, if we define polynomials $\mathcal{P}_n$ and $\mathcal{Q}_n$ by
\begin{equation*}
\frac{\td^n(\cot x)}{\td x^n}=\mathcal{P}_n(\cot x)\quad\text{and}\quad
\frac{\td^n(\csc x)}{\td x^n}=\mathcal{Q}_n(\cot x)\csc x,
\end{equation*}
then the polynomials $\mathcal{P}_n$ and $\mathcal{Q}_n$ for $n\in\mathbb{N}$ meet the recurrences
\begin{equation}
\begin{aligned}\label{recur-der-pol3}
\mathcal{P}_{n+1}(u)&=-\sum_{k=0}^n\binom{n}{k}\mathcal{P}_k(u)\mathcal{P}_{n-k}(u),\\
\mathcal{P}_{n+1}(u)&=-\bigl(1+u^2\bigr)\sum_{k=0}^n\binom{n}{k}\mathcal{Q}_k(u)\mathcal{Q}_{n-k}(u),\\
\mathcal{Q}_{n+1}(u)&=-\sum_{k=0}^n\binom{n}{k}\mathcal{P}_k(u)\mathcal{Q}_{n-k}(u).
\end{aligned}
\end{equation}
\par
What are the relations among $P_n$, $\mathcal{P}_n$, and $C_n$? What are the relations among $Q_n$, $\mathcal{Q}_n$, and $S_n$? Because they have similar recurrences in~\eqref{recur-der-polyn}, \eqref{recur-der-pol2}, and~\eqref{recur-der-pol3}, we conjecture that their differences are just a minus.
\par
Finally we just mention that the asymptotic distribution of zeros of some of the above derivative polynomials were investigated in~\cite{Dominici-VanAssche-AA2014}.

\section{Proofs of main results}

We are now in a position to prove Theorems~\ref{tan-deriv-thm} to~\ref{Ber-zeta-genoc-thm}.

\begin{proof}[Proof of Theorem~\ref{tan-deriv-thm}]
We prove this theorem by mathematical induction.
\par
It is easy to obtain that
\begin{equation*}
(\tan x)'=\sec^2x=\frac1{\cos^2x}\quad \text{and}\quad (\tan x)''=2\tan x\sec^2x=\frac{2\sin x}{\cos^3x}
\end{equation*}
This means $a_{1,0}=1$ and $a_{2,1}=2$.
Therefore, the formula~\eqref{tan-deriv-unity} is valid for $n=1,2$.
\par
Assume that the formulas~\eqref{tan-deriv-unity} and~\eqref{a(p-q)-unify} are valid for some $n>1$. By this inductive hypothesis and a direct differentiation, we have
\begin{align*}
\tan^{(2n+1)}x&=\bigl[\tan^{(2n)}x\bigr]'
=\Biggl\{\frac1{\cos^{2n+1}x}\sum_{\ell=0}^{n-1}a_{2n,2\ell+1}\sin[(2\ell+1)x]\Biggr\}'\\
&=\sum_{\ell=0}^{n-1}a_{2n,2\ell+1}\biggl\{\frac{\sin[(2\ell+1)x]}{\cos^{2n+1}x}\biggr\}'\\
&=\frac1{\cos^{2(n+1)}x}\sum_{\ell=0}^{n-1}a_{2n,2\ell+1}\{(n+1+\ell)\cos(2\ell x)+(\ell-n)\cos[2(\ell+1)x]\}\\
&=\frac1{\cos^{2(n+1)}x}\Biggl\{(n+1)a_{2n,1}+\sum_{\ell=1}^{n-1}[(n+1+\ell)a_{2n,2\ell+1}\\
&\quad-(n+1-\ell)a_{2n,2\ell-1}]\cos(2\ell x)-a_{2n,2n-1}\cos(2nx)\Biggr\}
\end{align*}
and
\begin{multline*}
\tan^{(2n+2)}x=\bigl[\tan^{(2n+1)}x\bigr]'
=\Biggl[\frac1{\cos^{2n+2}x}\sum_{\ell=0}^na_{2n+1,2\ell}\cos(2\ell x)\Biggr]'\\
\begin{aligned}
&=\sum_{\ell=0}^na_{2n+1,2\ell}\biggl[\frac{\cos(2\ell x)}{\cos^{2n+2}x}\biggr]'\\
&=\frac1{\cos^{2n+3}x}\sum_{\ell=0}^na_{2n+1,2\ell} \{(n+1-\ell)\sin[(2\ell+1)x]-(n+1+\ell)\sin[(2\ell-1)x]\}\\
&=\frac1{\cos^{2n+3}x}\Biggl\{[2(n+1)a_{2n+1,0}-(n+2)a_{2n+1,2}]\sin x+\sum_{\ell=1}^{n-1}[(n+1-\ell)a_{2n+1,2\ell}\\
&\quad-(n+2+\ell)a_{2n+1,2(\ell+1)}]\sin[(2\ell+1)x]+a_{2n+1,2n}\sin[(2n+1)x]\Biggr\}.
\end{aligned}
\end{multline*}
By straightforward computation, it is not difficult to see that
\begin{align*}
(n+1)a_{2n,1}&=a_{2n+1,0},\\
(n+1+\ell)a_{2n,2\ell+1}-(n+1-\ell)a_{2n,2\ell-1}&=a_{2n+1,2\ell},\\
-a_{2n,2n-1}&=a_{2n+1,2n},\\
2(n+1)a_{2n+1,0}-(n+2)a_{2n+1,2}&=a_{2n+2,1},\\
(n+1-\ell)a_{2n+1,2\ell}-(n+2+\ell)a_{2n+1,2(\ell+1)}&=a_{2n+2,2\ell+1},\\
a_{2n+1,2n}&=a_{2n+2,2n+1},
\end{align*}
where $1\le\ell\le n-1$. From these recurrence relations, the formulas~\eqref{tan-deriv-unity} and~\eqref{a(p-q)-unify} may be derived straightforwardly.
\par
The formula~\eqref{cotan-deriv-n} can be proved by induction as the proof of the formulas~\eqref{tan-deriv-unity} and~\eqref{a(p-q)-unify}. However, it is easy to be derived from the formulas~\eqref{tan-deriv-unity} and~\eqref{a(p-q)-unify} by considering $\cot x=-\tan\bigl(x+\frac\pi2\bigr)$ and $\cot^{(n)}x=-\tan^{(n)}\bigl(x+\frac\pi2\bigr)$ for $n\in\mathbb{N}$. Theorem~\ref{tan-deriv-thm} is thus proved.
\end{proof}

\begin{proof}[Proof of Theorem~\ref{bell-sin-cos-thm}]
In~\cite[p.~133]{Comtet-Combinatorics-74}, it was listed that
\begin{equation}\label{113-final-formula}
\frac1{k!}\Biggl(\sum_{m=1}^\infty x_m\frac{t^m}{m!}\Biggr)^k
=\sum_{n=k}^\infty \bell_{n,k}(x_1,x_2,\dotsc,x_{n-k+1})\frac{t^n}{n!}
\end{equation}
for $k\ge0$. Taking $x_m=\cos\bigl(x+m\frac{\pi}{2}\bigr)$ in~\eqref{113-final-formula} yields
\begin{equation}\label{bell-sin-cos-der}
\begin{gathered}
\sum_{n=k}^\infty \bell_{n,k}\biggl(-\sin x,-\cos x,\sin x,\cos x,\dotsc, -\sin\biggl[x+(n-k)\frac{\pi}{2}\biggr]\biggr)\frac{t^n}{n!}\\
=\frac1{k!}\Biggl[\sum_{m=1}^\infty \cos\biggl(x+m\frac{\pi}{2}\biggr)\frac{t^m}{m!}\Biggr]^k
=\frac1{k!}\biggl[-2 \sin \biggl(\frac{t}{2}\biggr) \sin \biggl(\frac{t}{2}+x\biggr)\biggr]^k\\
=\frac{[\cos(t+x)-\cos x]^k}{k!}
=\frac{(-1)^k}{k!}\sum_{\ell=0}^k(-1)^\ell\binom{k}{\ell}\cos^\ell(t+x)\cos^{k-\ell}x.
\end{gathered}
\end{equation}
Since $\cos x=\frac{e^{ix}+e^{-ix}}2$, we have
\begin{equation*}
(\cos x)^\ell=\frac1{2^\ell}(e^{ix}+e^{-ix})^\ell
=\frac1{2^\ell}\sum_{q=0}^\ell\binom{\ell}{q}e^{qix}e^{-(\ell-q)ix}
=\frac1{2^\ell}\sum_{q=0}^\ell\binom{\ell}{q}e^{(2q-\ell)ix}
\end{equation*}
and
\begin{gather*}
\frac{\td^m(\cos x)^\ell}{\td x^m}
=\frac1{2^\ell}\sum_{q=0}^\ell\binom{\ell}{q}(2q-\ell)^mi^me^{(2q-\ell)ix}
=\frac1{2^\ell}\sum_{q=0}^\ell\binom{\ell}{q}(2q-\ell)^me^{[\pi m/2+(2q-\ell)x]i}\\
=\frac1{2^\ell}\sum_{q=0}^\ell\binom{\ell}{q}(2q-\ell)^m \biggl\{\cos\biggl[\frac{\pi}2m+(2q-\ell)x\biggr] +i\sin\biggl[\frac{\pi}2m+(2q-\ell)x\biggr]\biggr\}
\end{gather*}
which implies the formula~\eqref{cos-ell-derr} and the equation~\eqref{eq=0-1-sin-cos}.
\par
Differentiating $m\ge k$ times with respect to $t$ on the very ends of~\eqref{bell-sin-cos-der} and employing~\eqref{cos-ell-derr} result in
\begin{align*}
&\quad\sum_{n=m}^\infty \bell_{n,k}\biggl(-\sin x,-\cos x,\sin x,\cos x,\dotsc, -\sin\biggl[x+(n-k)\frac{\pi}{2}\biggr]\biggr)\frac{t^{n-m}}{(n-m)!}\\
&=\frac{(-1)^k}{k!}\sum_{\ell=0}^k(-1)^\ell\binom{k}{\ell}
\frac{\td^m\cos^\ell(t+x)}{\td t^m}\cos^{k-\ell}x\\
&=\frac{(-1)^k}{k!}\sum_{\ell=0}^k(-1)^\ell\binom{k}{\ell}
\frac1{2^\ell}\sum_{q=0}^\ell\binom{\ell}{q}(2q-\ell)^m \cos\biggl[\frac{\pi}2m+(2q-\ell)(t+x)\biggr]\cos^{k-\ell}x.
\end{align*}
Further taking the limit $t\to0$ and rearranging reveal~\eqref{bell-sin-eq}.
\par
Taking $x_m=\sin\bigl(x+m\frac{\pi}{2}\bigr)$ in~\eqref{113-final-formula} yields
\begin{equation}\label{bell-sin-der}
\begin{gathered}
\sum_{n=k}^\infty \bell_{n,k}\biggl(\cos x,-\sin x,-\cos x,\sin x,\dotsc, -\cos\biggl[x+(n-k)\frac{\pi}{2}\biggr]\biggr)\frac{t^n}{n!}\\
=\frac1{k!}\Biggl[\sum_{m=1}^\infty \sin\biggl(x+m\frac{\pi}{2}\biggr)\frac{t^m}{m!}\Biggr]^k
=\frac1{k!}\biggl[2 \sin \biggl(\frac{t}{2}\biggr) \cos\biggl(\frac{t}{2}+x\biggr)\biggr]^k\\
=\frac{[\sin(t+x)-\sin x]^k}{k!}
=\frac{(-1)^k}{k!}\sum_{\ell=0}^k(-1)^\ell\binom{k}{\ell}\sin^\ell(t+x)\sin^{k-\ell}x.
\end{gathered}
\end{equation}
Since $\sin x=\frac{e^{ix}-e^{-ix}}{2i}$, we have
\begin{equation*}
(\sin x)^\ell
=\frac1{(2i)^\ell}\sum_{q=0}^\ell(-1)^{\ell-q}\binom{\ell}{q}e^{qix}e^{-(\ell-q)ix}
=\frac{(-1)^\ell}{(2i)^\ell}\sum_{q=0}^\ell(-1)^q\binom{\ell}{q}e^{(2q-\ell)ix}
\end{equation*}
and
\begin{gather*}
\begin{aligned}
\frac{\td^m(\sin x)^\ell}{\td x^m}
&=\frac{(-1)^\ell}{(2i)^\ell} \sum_{q=0}^\ell(-1)^q\binom{\ell}{q}(2q-\ell)^mi^me^{(2q-\ell)ix}\\
&=\frac{(-1)^\ell}{2^\ell} \sum_{q=0}^\ell(-1)^q\binom{\ell}{q}(2q-\ell)^mi^{m-\ell} e^{(2q-\ell)ix}\\
&=\frac{(-1)^\ell}{2^\ell} \sum_{q=0}^\ell(-1)^q\binom{\ell}{q}(2q-\ell)^m e^{[(m-\ell)\pi/2+(2q-\ell)x]i}
\end{aligned}\\
=\frac{(-1)^\ell}{2^\ell} \sum_{q=0}^\ell(-1)^q\binom{\ell}{q}(2q-\ell)^m \biggl\{\cos\biggl[(m-\ell)\frac\pi2+(2q-\ell)x\biggr] +i\sin\biggl[(m-\ell)\frac\pi2+(2q-\ell)x\biggr]\biggr\}
\end{gather*}
which means the formula~\eqref{sin-ell-d} and the equation~\eqref{eq=0-2-sin-cos}.
\par
Differentiating $m\ge k$ times with respect to $t$ on the very ends of~\eqref{bell-sin-der} and utilizing~\eqref{sin-ell-d} lead to
\begin{multline*}
\sum_{n=m}^\infty \bell_{n,k}\biggl(\cos x,-\sin x,-\cos x,\sin x,\dotsc, -\cos\biggl[x+(n-k)\frac{\pi}{2}\biggr]\biggr)\frac{t^{n-m}}{(n-m)!}\\
\begin{aligned}
&=\frac{(-1)^k}{k!}\sum_{\ell=0}^k(-1)^\ell\binom{k}{\ell}
\frac{\td^m\sin^\ell(t+x)}{\td x^m}\sin^{k-\ell}x\\
&=\frac{(-1)^k}{k!}\sum_{\ell=0}^k\binom{k}{\ell}
\frac1{2^\ell} \sum_{q=0}^\ell(-1)^q\binom{\ell}{q}(2q-\ell)^m \cos\biggl[(m-\ell)\frac\pi2+(2q-\ell)(t+x)\biggr]\sin^{k-\ell}x.
\end{aligned}
\end{multline*}
Further letting $t\to0$ gives the formula~\eqref{bell-sin=ans}.
\par
Applying $f(u)=\ln u$ and $u=h(t)=\cos t$ in~\eqref{Bruno-Bell-Polynomial} and using the formula~\eqref{bell-sin-eq} give
\begin{multline*}
\frac{\td^n\tan t}{\td t^n}=-\frac{\td^{n+1}\ln(\cos t)}{\td t^{n+1}}
=-\sum_{k=1}^{n+1}(\ln u)^{(k)}\bell_{n+1,k} \bigl(\cos't,\cos''t,\dotsc,\cos^{(n-k+2)}t\bigr)\\
=-\sum_{k=1}^{n+1}\frac{(-1)^{k-1}(k-1)!}{u^k} \bell_{n+1,k}\biggl(-\sin t,-\cos t,\sin t,\cos t,\dotsc, -\cos\biggl[t+(n-k)\frac{\pi}{2}\biggr]\biggr)\\
=\sum_{k=1}^{n+1}\frac1{k}\sum_{\ell=0}^k\frac{(-1)^\ell}{2^\ell}\binom{k}{\ell}
\frac1{\cos^{\ell}t}\sum_{q=0}^\ell\binom{\ell}{q}(2q-\ell)^{n+1} \cos\biggl[\frac{\pi}2(n+1)+(2q-\ell)t\biggr].
\end{multline*}
The formula~\eqref{tan-der-nth-eq} follows.
\par
Applying $f(u)=\ln u$ and $u=h(t)=\sin t$ in~\eqref{Bruno-Bell-Polynomial} and using the formula~\eqref{bell-sin=ans} generate
\begin{multline*}
\frac{\td^n\cot t}{\td t^n}=\frac{\td^{n+1}\ln(\sin t)}{\td t^{n+1}}
=\sum_{k=1}^{n+1}(\ln u)^{(k)}\bell_{n+1,k} \bigl(\sin't,\sin''t,\dotsc,\sin^{(n-k+2)}t\bigr)\\
=\sum_{k=1}^{n+1}\frac{(-1)^{k-1}(k-1)!}{u^k} \bell_{n+1,k}\biggl(\cos t,-\sin t,-\cos t,\sin t,\dotsc, -\sin\biggl[t+(n-k)\frac{\pi}{2}\biggr]\biggr)\\
=\sum_{k=1}^{n+1}\frac1k\sum_{\ell=0}^k\frac1{2^\ell}\binom{k}{\ell}\frac1{\sin^{\ell}t}
\sum_{q=0}^\ell(-1)^q\binom{\ell}{q}(2q-\ell)^{n+1} \sin\biggl[(n-\ell)\frac\pi2+(2q-\ell)t\biggr].
\end{multline*}
The formula~\eqref{cotan-der-nth-eq} follows.
The proof of Theorem~\ref{bell-sin-cos-thm} is complete.
\end{proof}

\begin{proof}[Proof of Theorem~\ref{tan-cot-app-Bernoulli}]
Taking $n=2m-1$ and letting $x\to0$ in~\eqref{tan-deriv-unity} yield
\begin{equation*}
T_{2m-1}=\lim_{x\to0}\tan^{(2m-1)}x
=\frac12a_{2m-1,0}
+\sum_{k=1}^{m-1} a_{2m-1,2k},
\end{equation*}
where
\begin{equation*}
a_{2m-1,0}=2\sum_{\ell=0}^{m-1} (-1)^{\ell}\binom{2m}{\ell}(m-\ell)^{2m-1}
\end{equation*}
and
\begin{equation*}
a_{2m-1,2k}=(-1)^k2\sum_{\ell=0}^{m-k-1} (-1)^{\ell}\binom{2m}{\ell}(m-k-\ell)^{2m-1}.
\end{equation*}
The formula~\eqref{Tangent-No-eq} follows.
\par
Taking $x\to0$ in~\eqref{tan-der-nth-eq} gives
\begin{align*}
\lim_{x\to0}\frac{\td^{2m-1}\tan x}{\td x^{2m-1}}
&=-\sum_{k=1}^{2m}\frac1{k}\sum_{\ell=0}^k\frac{(-1)^\ell}{2^\ell}\binom{k}{\ell}
\sum_{q=0}^\ell\binom{\ell}{q}(2q-\ell)^{2m} \sin\biggl[\frac{\pi}2(2m-1)\biggr]\\
&=(-1)^m\sum_{k=1}^{2m}\frac1{k}\sum_{\ell=0}^k(-1)^\ell\frac1{2^\ell}\binom{k}{\ell}
\sum_{q=0}^\ell\binom{\ell}{q}(2q-\ell)^{2m}\\
&=(-1)^m\sum_{\ell=1}^{2m}(-1)^\ell\frac1{2^\ell} \sum_{q=0}^\ell\binom{\ell}{q}(2q-\ell)^{2m} \sum_{k=\ell}^{2m}\frac1{k}\binom{k}{\ell}\\
&=(-1)^m\sum_{\ell=1}^{2m}(-1)^\ell\frac1{2^\ell} \frac{2m+1-\ell}{(2m+1)\ell}\binom{2m+1}{\ell} \sum_{q=0}^\ell\binom{\ell}{q}(2q-\ell)^{2m}\\
&=T_{2m-1}.
\end{align*}
The formula~\eqref{tan-No-faaa} follows.
\par
Taking $x\to0$ in~\eqref{tan-deriv-Quadratic-form} figures out
\begin{gather*}
T_{2m+1}=\lim_{x\to0}y^{(2m+1)}
=\sum_{k=1}^m\binom{2m}{2k-1}\lim_{x\to0}\bigl[y^{(2k-1)}y^{(2m-2k+1)}\bigr]
=\sum_{k=1}^m\binom{2m}{2k-1}T_{2k-1}T_{2(m-k)+1}.
\end{gather*}
The recurrence~\eqref{tan-No-rec-eq} follows.
The proof of Theorem~\ref{tan-cot-app-Bernoulli} is thus proved.
\end{proof}

\begin{proof}[Proof of Theorem~\ref{Ber-zeta-genoc-thm}]
From~\eqref{tan-Bern-lim-der}, it follows that
\begin{equation*}
B_{2k}=\frac{(-1)^{k-1}k} {2^{2k-1}\bigl(2^{2k}-1\bigr)}T_{2k-1}.
\end{equation*}
Combining this with the formulas~\eqref{Tangent-No-eq} and~\eqref{tan-No-faaa} gains the formulas~\eqref{Tangent-No-eq-B} and~\eqref{tan-No-faaa-B}.
\par
Substituting~\eqref{tan-Bern-lim-der} into~\eqref{tan-No-rec-eq} and rearranging yield~\eqref{tan-No-rec-eq-B}.
\par
In~\cite[p.~807, 23.2.16]{abram}, it is listed that for $n\ge1$
\begin{equation}\label{zeta(2q)-B(2q)}
\zeta(2n)=\frac{(2\pi)^{2n}}{2(2n)!}|B_{2n}|.
\end{equation}
See also~\cite[p.~332]{Remmert-GTM}. Substituting the formulas~\eqref{Tangent-No-eq-B} and~\eqref{tan-No-faaa-B} into~\eqref{zeta(2q)-B(2q)} results in~\eqref{zeta-No-eq-a} and~\eqref{zeta-No-eq-b}.
\par
Substituting~\eqref{zeta(2q)-B(2q)} into~\eqref{tan-No-rec-eq-B} and rearranging yield~\eqref{zeta-No-rec-eq-B}.
\par
The Genocchi numbers satisfy $G_1=1$, $G_{2n+1}=0$, and
\begin{equation*}
G_{2n}=2\bigl(1-2^{2n}\bigr)B_{2n}=2nE_{2n-1}(0)
\end{equation*}
for $n\in\mathbb{N}$. By similar arguments as above, we may obtain the formulas~\eqref{G-No-eq-a}, \eqref{G-No-eq-b}, \eqref{E-No-eq-a}, and~\eqref{E-No-eq-b} and the recurrence relations~\eqref{No-rec-eq-G} and~\eqref{tan-No-rec-eq-E}.
The proof of Theorem~\ref{Ber-zeta-genoc-thm} is complete.
\end{proof}

\section{Remarks}

Finally, we give several remarks on related results.

\begin{remark}
Let $y=\tan x$. The relation~\eqref{tan-deriv-Quadratic-form} may be rewritten as
\begin{equation*}
y^{(n)}=YAY^T
\end{equation*}
for $n\ge2$, where $Y^T$ stands for the transpose of the $1\times n$ matrix
\begin{equation*}
Y=\bigl(y,y',y'',\dotsc,y^{(n-2)}, y^{(n-1)}\bigr)
\end{equation*}
and $A$ is a $n\times n$ square matrix whose elements $a_{k,\ell}$ satisfy
\begin{equation*}
a_{k,\ell}=
\begin{cases}
\dbinom{n-1}{\ell-1}, & k+\ell=n+1,\\
0, & k+\ell\ne n+1.
\end{cases}
\end{equation*}
\end{remark}

\begin{remark}
By the formulas in~\eqref{a-b-relation-ratio}, the tangent numbers $T_{2m-1}$ for $m\in\mathbb{N}$ may be represented as
\begin{equation*}
T_{2m-1}=
-\begin{vmatrix}
0 & 1 & 0 & 0 & \dots & 0\\
1 & 0 & 1 & 0 & \dots & 0\\
0 & -1 & 0 & 1 & \dots & 0\\
-1 & 0 & -\binom{3}{1} & 0 & \dots & 0\\
\hdotsfor[2]{6}\\
0 & (-1)^{m-1} & 0 & (-1)^{m-2}\binom{2m-2}2 & \dots & 1\\
(-1)^{m-1} & 0 & (-1)^{m-1}\binom{2m-1}1 & 0 & \dots & 0
\end{vmatrix}.
\end{equation*}
\end{remark}

\begin{remark}
The general and explicit formula~\eqref{cotan-deriv-n} for derivatives of the cotangent function $\cot x$ in Theorem~\ref{tan-deriv-thm} has been applied in~\cite{singularity-combined.tex} to establish limit formulas for ratios of polygamma functions at their singularities.
\end{remark}

\begin{remark}
The formulas~\eqref{zeta-No-eq-a}, \eqref{zeta-No-eq-b}, and~\eqref{zeta-No-rec-eq-B} show the irrationality of $\zeta(2n)$ for $n\in\mathbb{N}$. However, a problem about the irrationality of $\zeta(2n+1)$ for $n\ge2$ is still keeping open, to the best of the author's knowledge. For more information, please refer to~\cite{Zeta-luo.tex, zeta(3)} and closely related references therein.
\end{remark}

\begin{remark}
In~\cite{Wallis-Ratio-Sum.tex}, the method and idea to prove the formulas~\eqref{sin-ell-d} and~\eqref{cos-ell-derr} were employed to express the Wallis ratio formula as a sum.
\end{remark}

\begin{remark}
In~\cite{exp-derivative-sum-Combined.tex}, among other things, by establishing the derivative polynomials~\eqref{der-polynomial-exp} for $\lambda=\alpha=1$, that is,
\begin{equation*}
\biggl(\frac1{e^t-1}\biggr)^{(q)}=(-1)^q\sum_{m=1}^{q+1}(m-1)!S(q+1,m)\biggl(\frac1{e^t-1}\biggr)^m
\end{equation*}
for $q\in\{0\}\cup\mathbb{N}$, an explicit formula
\begin{equation*}
B_{2k}=1+\sum_{m=1}^{2k-1}\frac{S(2k+1,m+1) S(2k,2k-m)}{\binom{2k}{m}}
-\frac{2k}{2k+1}\sum_{m=1}^{2k}\frac{S(2k,m)S(2k+1,2k-m+1)}{\binom{2k}{m-1}},\quad k\in\mathbb{N}
\end{equation*}
for computing the Bernoulli numbers $B_{2k}$ in terms of the Stirling numbers $S(q,m)$ of the second kind was derived. For more information on closely related problems, please refer to~\cite{Eight-Identy-More.tex-JCAM, CAM-D-13-01430-Xu-Cen} and the references cited therein.
\end{remark}

\begin{remark}
In~\cite{Filomat-36-73-1.tex}, based on establishment of the general and explicit formula
\begin{equation*}
\biggl(\frac1{\ln x}\biggr)^{(n)}=\frac{(-1)^n}{x^n}\sum_{q=2}^{n+1}\frac{\alpha_{n,q}}{(\ln x)^{q}}, \quad n\in\mathbb{N},
\end{equation*}
where
\begin{equation}\label{a=n=0}
\alpha_{n,2}=(n-1)!
\end{equation}
and, for $n+1\ge q\ge3$,
\begin{equation}\label{a=n=q=eq}
\alpha_{n,q}=(q-1)!(n-1)!\sum_{\ell_1=1}^{n-1} \frac1{\ell_1}\sum_{\ell_2=1}^{\ell_1-1}\frac1{\ell_2}\dotsm \sum_{\ell_{q-3}=1}^{\ell_{q-4}-1}\frac1{\ell_{q-3}} \sum_{\ell_{q-2}=1}^{\ell_{q-3}-1}\frac1{\ell_{q-2}},
\end{equation}
among other things, it was discovered that the Stirling numbers of the first kind $s(n,q)$ for $1\le q\le n$ may be computed by
\begin{equation*}
s(n,q)=(-1)^{n+q}(n-1)!\sum_{\ell_1=1}^{n-1} \frac1{\ell_1}\sum_{\ell_2=1}^{\ell_1-1}\frac1{\ell_2}\dotsm \sum_{\ell_{q-2}=1}^{\ell_{q-3}-1}\frac1{\ell_{q-2}} \sum_{\ell_{q-1}=1}^{\ell_{q-2}-1}\frac1{\ell_{q-1}},
\end{equation*}
equivalently,
\begin{equation*}
(-1)^{n-k}\frac{s(n,k)}{(n-1)!}= \sum _{m=k-1}^{n-1}\frac1m\biggl[(-1)^{m-(k-1)}\frac{s(m,k-1)}{(m-1)!}\biggr],\quad 2\le k\le n,
\end{equation*}
and that the Bernoulli numbers of the second kind $b_n$ for $n\ge2$ may be computed by
\begin{equation*}
b_n=(-1)^n\frac1{n!}\Biggl(\frac1{n+1}+\sum_{k=2}^n \frac{\alpha_{n,k}-n\alpha_{n-1,k}}{k!}\Biggr),
\end{equation*}
where the Stirling numbers of the first kind $s(n,k)$ and the Bernoulli numbers of the second kind $b_n$ for $n\ge0$ may be generated respectively by
\begin{equation*}
\frac{[\ln(1+x)]^m}{m!}=\sum_{k=m}^\infty\frac{s(k,m)}{k!}x^k,\quad |x|<1
\end{equation*}
and
\begin{equation*}
\frac{x}{\ln(1+x)}=\sum_{n=0}^\infty b_nx^n,\quad |x|<1,
\end{equation*}
and $\alpha_{n,k}$ are defined by~\eqref{a=n=0} and~\eqref{a=n=q=eq}.
\end{remark}

\begin{remark}
In~\cite{1st-Sirling-Number-2012.tex}, three integral representations for the Stirling numbers of the first kind $s(n,k)$ were created, one of them is
\begin{equation}\label{n-times-diriv}
s(n,k)=\binom{n}{k}\lim_{x\to0}\frac{\td^{n-k}}{\td x^{n-k}} \biggl\{\biggl[\int_0^\infty\biggl(\int_{1/e}^1 t^{xu-1}\td t\biggr)e^{-u}\td u\biggr]^k\biggr\}
\end{equation}
for $1\le k\le n$. By virtue of~\eqref{n-times-diriv}, the sequence
\begin{equation*}
\Biggl\{(-1)^n\frac{s(n+k,k)}{\binom{n+k}{k}}\Biggr\}_{n\ge0}
\end{equation*}
for any given $k\in\mathbb{N}$ is proved to be logarithmically convex.
\end{remark}

\begin{remark}
In~\cite{KMS-B14-0477.tex}, it was obtained that the Bernoulli numbers of the second kind $b_n$ may be represented by
\begin{equation}\label{Bernoulli-2rd-int}
b_n=(-1)^{n+1}\int_1^\infty\frac{1}{\{[\ln(t-1)]^2+\pi^2\}t^n}\td t, \quad n\in\mathbb{N}.
\end{equation}
As a result of the integral representation~\eqref{Bernoulli-2rd-int}, the sequence $\{(-1)^nb_{n+1}\}_{n\ge0}$ of the Bernoulli numbers of the second kind was proved to be completely monotonic. In~\cite{KMS-B14-0477.tex}, among other things, the sequence $\{(-1)^qq!b_{q+1}\}_{q\ge0}$ is proved to be logarithmically convex.
\end{remark}

\begin{remark}
The general and explicit formulas for the $n$-th derivatives of the exponential function $e^{\pm1/x}$ were established in~\cite{DMST-MM-2013-Exp, exp-reciprocal-cm-IJOPCM.tex} and applied in~\cite{Bell-Stirling-Lah-simp.tex, Filomat-36-73-1.tex, Bessel-ineq-Dgree-CM.tex, QiBerg.tex, simp-exp-degree-revised.tex} and related references therein.
\end{remark}

\begin{remark}
Recently, some results on explicit formulas and inequalities for the Bell numbers, the Bell polynomials $\bell_{n,k}$, the Bernoulli numbers $B_{2n}$, the Cauchy numbers, the Euler polynomials $E_{n}$, and the Stirling numbers were discovered or recovered in~\cite{Genocchi-Stirling.tex, ANLY-D-12-1238.tex, Bell-Stirling-HyperGeom.tex, Guo-Qi-JANT-Bernoulli.tex, Eight-Identy-More.tex-JCAM, Special-Bell2Euler.tex, exp-derivative-sum-Combined.tex, Lah-Num-Int-Prop-Corr.tex, Bernoulli-Ratio-Ineq.tex-Gate, Bernoulli-Ratio-Ineq.tex-RGMIA, Norlund-No-CM-JNT.tex, 1st-Sirling-Number-2012.tex, Bernoulli-Stirling-4P.tex} and closely related references therein.
\end{remark}

\begin{remark}
This paper is a corrected, expanded, and revised version of the preprint~\cite{derivative-tan-cot.tex}.
\end{remark}

\subsection*{Acknowledgments}
The author thanks Professor Sergei M. Sitnik in Voronezh Institute of the Ministry of Internal Affairs of Russia for his providing the formula~\eqref{Sitnik-Bourbaki} and the reference~\cite{Bourbaki-Spain-2004} on September~25, 2014.
\par
The author appreciates anonymous referees for their careful corrections to and valuable comments on the original version of this paper.


\begin{thebibliography}{99}

\bibitem{abram}
M. Abramowitz and I. A. Stegun (Eds), \textit{Handbook of Mathematical Functions with Formulas, Graphs, and Mathematical Tables}, National Bureau of Standards, Applied Mathematics Series \textbf{55}, 10th printing, Washington, 1972.

\bibitem{Bourbaki-Spain-2004}
N. Bourbaki, \emph{Functions of a Real Variable, Elementary Theory}, Translated from the 1976 French original by Philip Spain. Elements of Mathematics (Berlin). Springer-Verlag, Berlin, 2004; Available online at \url{http://dx.doi.org/10.1007/978-3-642-59315-4}.

\bibitem{Boyadzhiev-Fibonacci2007}
K. N. Boyadzhiev, \emph{Derivative polynomials for tanh, tan, sech and sec in explicit form}, Fibonacci Quart. \textbf{45} (2007), no.~4, 291\nobreakdash--303.

\bibitem{Boyadzhiev-arXiv0903}
K. N. Boyadzhiev, \emph{Derivative Polynomials for tanh, tan, sech and sec in explicit form}, available online at \url{http://arxiv.org/abs/0903.0117}.

\bibitem{Comtet-Combinatorics-74}
L. Comtet, \emph{Advanced Combinatorics: The Art of Finite and Infinite Expansions}, Revised and Enlarged Edition, D. Reidel Publishing Co., 1974.

\bibitem{DMST-MM-2013-Exp}
S. Daboul, J. Mangaldan, M. Z. Spivey, and P. J. Taylor, \emph{The Lah numbers and the $n$th derivative of $e\sp {1/x}$}, Math. Mag. \textbf{86} (2013), no.~1, 39\nobreakdash--47; Available online at \url{http://dx.doi.org/10.4169/math.mag.86.1.039}.

\bibitem{Dominici-VanAssche-AA2014}
D. Dominici and W. Van Assche, \emph{Zero distribution of polynomials satisfying a differential-difference equation}, Anal. Appl. (Singap.) \textbf{12} (2014), no.~6, 635\nobreakdash--666; Available online at \url{http://dx.doi.org/10.1142/S0219530514500390}.

\bibitem{Genocchi-Stirling.tex}
B.-N. Guo and F. Qi, \textit{A new explicit formula for the Bernoulli and Genocchi numbers in terms of the Stirling numbers}, Glob. J. Math. Anal. \textbf{3} (2015), no.~1, 33\nobreakdash--36; Available online at \url{http://dx.doi.org/10.14419/gjma.v3i1.4168}.

\bibitem{ANLY-D-12-1238.tex}
B.-N. Guo and F. Qi, \textit{Alternative proofs of a formula for Bernoulli numbers in terms of Stirling numbers}, Analysis (Berlin) \textbf{34} (2014), no.~2, 187\nobreakdash--193; Available online at \url{http://dx.doi.org/10.1515/anly-2012-1238}.

\bibitem{Bell-Stirling-HyperGeom.tex}
B.-N. Guo and F. Qi, \textit{An explicit formula for Bell numbers in terms of Stirling numbers and hypergeometric functions}, Glob. J. Math. Anal. \textbf{2} (2014), no.~4, 243\nobreakdash--248; Available online at \url{http://dx.doi.org/10.14419/gjma.v2i4.3310}.

\bibitem{Guo-Qi-JANT-Bernoulli.tex}
B.-N. Guo and F. Qi, \textit{An explicit formula for Bernoulli numbers in terms of Stirling numbers of the second kind}, J. Anal. Number Theory \textbf{3} (2015), no.~1, 27\nobreakdash--30; Available online at \url{http://dx.doi.org/10.12785/jant/030105}.

\bibitem{Eight-Identy-More.tex-JCAM}
B.-N. Guo and F. Qi, \textit{Explicit formulae for computing Euler polynomials in terms of Stirling numbers of the second kind}, J. Comput. Appl. Math. \textbf{272} (2014), 251\nobreakdash--257; Available online at \url{http://dx.doi.org/10.1016/j.cam.2014.05.018}.

\bibitem{Special-Bell2Euler.tex}
B.-N. Guo and F. Qi, \textit{Explicit formulas for special values of the Bell polynomials of the second kind and the Euler numbers}, ResearchGate Technical Report, available online at \url{http://dx.doi.org/10.13140/2.1.3794.8808}.

\bibitem{Wallis-Ratio-Sum.tex}
B.-N. Guo and F. Qi, \textit{On the Wallis formula}, Internat. J. Anal. Appl. \textbf{8} (2015), no.~1, 30\nobreakdash--38.

\bibitem{exp-derivative-sum-Combined.tex}
B.-N. Guo and F. Qi, \textit{Some identities and an explicit formula for Bernoulli and Stirling numbers}, J. Comput. Appl. Math. \textbf{255} (2014), 568\nobreakdash--579; Available online at \url{http://dx.doi.org/10.1016/j.cam.2013.06.020}.

\bibitem{Lah-Num-Int-Prop-Corr.tex}
B.-N. Guo and F. Qi, \textit{Some integral representations and properties of Lah numbers}, J. Algebra Number Theory Acad. \textbf{4} (2014), no.~3, 77\nobreakdash--87.

\bibitem{Hoffman-Monthly1995}
M. E. Hoffman, \emph{Derivative polynomials for tangent and secant}, Amer. Math. Monthly \textbf{102} (1995), no.~1, 23\nobreakdash--30; Available online at \url{http://dx.doi.org/10.2307/2974853}.

\bibitem{Hoffman-EJC1999}
M. E. Hoffman, \emph{Derivative polynomials, Euler polynomials, and associated integer sequences}, Electron. J. Combin. \textbf{6} (1999), Research Paper~21, 13~pages; Available online at \url{http://www.combinatorics.org/Volume_6/Abstracts/v6i1r21.html}.

\bibitem{Zeta-luo.tex}
Q.-M. Luo, B.-N. Guo, and F. Qi, \textit{On evaluation of Riemann zeta function $\zeta(s)$}, Adv. Stud. Contemp. Math. (Kyungshang) \textbf{7} (2003), no.~2, 135\nobreakdash--144.

\bibitem{zeta(3)}
Q.-M. Luo, Z.-L. Wei, and F. Qi, \textit{Lower and upper bounds of $\zeta(3)$}, Adv. Stud. Contemp. Math. (Kyungshang) \textbf{6} (2003), no.~1, 47\nobreakdash--51.

\bibitem{Ma-EJC-2013}
S.-M. Ma, \emph{A family of two-variable derivative polynomials for tangent and secant}, Electron. J. Combin. \textbf{20} (2013), no.~1, Paper~11, 12 pages.

\bibitem{Ma-BAMS-2014}
S.-M. Ma, \emph{On $\gamma$-vectors and the derivatives of the tangent and secant functions}, Bull. Aust. Math. Soc. \textbf{90} (2014), no.~2, 177\nobreakdash--185; Available online at \url{http://dx.doi.org/10.1017/S0004972714000057}.

\bibitem{MA-EJC}
S.-M. Ma, \emph{Some combinatorial arrays generated by context-free grammars}, European J. Combin. \textbf{34} (2013), no.~7, 1081\nobreakdash--1091; Available online at \url{http://dx.doi.org/10.1016/j.ejc.2013.03.002}.

\bibitem{Ma-Mansour-Wang-1403}
S.-M. Ma, T. Mansour, and D. G. L. Wang, \textit{Combinatorics of Dumont differential system on the  Jacobi elliptic functions}, available online at \url{http://arxiv.org/abs/1403.0233}.

\bibitem{Bernoulli-Ratio-Ineq.tex-Gate}
F. Qi, \textit{A double inequality for ratios of Bernoulli numbers}, ResearchGate Dataset, available online at \url{http://dx.doi.org/10.13140/2.1.2367.2962}.

\bibitem{Bernoulli-Ratio-Ineq.tex-RGMIA}
F. Qi, \textit{A double inequality for ratios of Bernoulli numbers}, RGMIA Res. Rep. Coll. \textbf{17} (2014), Article~103, 4~pages; Available online at \url{http://rgmia.org/v17.php}.

\bibitem{Bell-Stirling-Lah-simp.tex}
F. Qi, \textit{An explicit formula for computing Bell numbers in terms of Lah and Stirling numbers}, available online at \url{http://arxiv.org/abs/1401.1625}.

\bibitem{Norlund-No-CM-JNT.tex}
F. Qi, \textit{An integral representation, complete monotonicity, and inequalities of Cauchy numbers of the second kind}, J. Number Theory \textbf{144} (2014), 244\nobreakdash--255; Available online at \url{http://dx.doi.org/10.1016/j.jnt.2014.05.009}.

\bibitem{Filomat-36-73-1.tex}
F. Qi, \textit{Explicit formulas for computing Bernoulli numbers of the second kind and Stirling numbers of the first kind}, Filomat \textbf{28} (2014), no.~2, 319\nobreakdash--327; Available online at \url{http://dx.doi.org/10.2298/FIL1402319O}.

\bibitem{derivative-tan-cot.tex}
F. Qi, \textit{Explicit formulas for the $n$-th derivatives of the tangent and cotangent functions}, available online at \url{http://arxiv.org/abs/1202.1205}.

\bibitem{1st-Sirling-Number-2012.tex}
F. Qi, \textit{Integral representations and properties of Stirling numbers of the first kind}, J. Number Theory \textbf{133} (2013), no.~7, 2307\nobreakdash--2319; Available online at \url{http://dx.doi.org/10.1016/j.jnt.2012.12.015}.

\bibitem{singularity-combined.tex}
F. Qi, \textit{Limit formulas for ratios between derivatives of the gamma and digamma functions at their singularities}, Filomat \textbf{27} (2013), no.~4, 601\nobreakdash--604; Available online at \url{http://dx.doi.org/10.2298/FIL1304601Q}.

\bibitem{Bessel-ineq-Dgree-CM.tex}
F. Qi, \textit{Properties of modified Bessel functions and completely monotonic degrees of differences between exponential and trigamma functions}, Math. Inequal. Appl. \textbf{18} (2015), no.~2, 493\nobreakdash--518; Available online at \url{http://dx.doi.org/10.7153/mia-18-37}.

\bibitem{QiBerg.tex}
F. Qi and C. Berg, \textit{Complete monotonicity of a difference between the exponential and trigamma functions and properties related to a modified Bessel function}, Mediterr. J. Math. \textbf{10} (2013), no.~4, 1685\nobreakdash--1696; Available online at \url{http://dx.doi.org/10.1007/s00009-013-0272-2}.

\bibitem{Bernoulli-Stirling-4P.tex}
F. Qi and B.-N. Guo, \textit{Alternative proofs of a formula for Bernoulli numbers in terms of Stirling numbers}, Analysis (Berlin) \textbf{34} (2014), no.~3, 311\nobreakdash--317; Available online at \url{http://dx.doi.org/10.1515/anly-2014-0003}.

\bibitem{simp-exp-degree-revised.tex}
F. Qi and S.-H. Wang, \textit{Complete monotonicity, completely monotonic degree, integral representations, and an inequality related to the exponential, trigamma, and modified Bessel functions}, Glob. J. Math. Anal. \textbf{2} (2014), no.~3, 91\nobreakdash--97; Available online at \url{http://dx.doi.org/10.14419/gjma.v2i3.2919}.

\bibitem{Euler-No-3Sum.tex}
F. Qi and C.-F. Wei, \textit{Several closed expressions for the Euler numbers}, ResearchGate Technical Report, available online at \url{http://dx.doi.org/10.13140/2.1.3474.7688}.

\bibitem{KMS-B14-0477.tex}
F. Qi and X.-J. Zhang, \textit{An integral representation, some inequalities, and complete monotonicity of the Bernoulli numbers of the second kind}, Bull. Korean Math. Soc. \textbf{52} (2015), no.~3, 987\nobreakdash--998; Available online at \url{http://dx.doi.org/10.4134/BKMS.2015.52.3.987}.

\bibitem{Deriv-Arcs-Cos.tex}
F. Qi and M.-M. Zheng, \textit{Explicit expressions for a family of the Bell polynomials and applications}, Appl. Math. Comput. \textbf{258} (2015), 597\nobreakdash--607; Available online at \url{http://dx.doi.org/10.1016/j.amc.2015.02.027}.

\bibitem{Remmert-GTM}
R. Remmert, \textit{Theory of Complex Functions}, Translated from the second German edition by Robert B. Burckel, Graduate Texts in Mathematics 122, Readings in Mathematics, Springer-Verlag, New York, 1991; Available online at \url{http://dx.doi.org/10.1007/978-1-4612-0939-3}.

\bibitem{Schwatt-1924}
I. J. Schwatt, \emph{An Introduction to the Operations with Series}, Chelsea Publishing Co., New York, 1924; Available online at \url{http://hdl.handle.net/2027/wu.89043168475}.

\bibitem{Schwatt-1962}
I. J. Schwatt, \emph{An Introduction to the Operations with Series}, Second edition, Chelsea Publishing Co., New York, 1962.

\bibitem{CAM-D-13-01430-Xu-Cen}
A.-M. Xu and Z.-D. Cen, \textit{Some identities involving exponential functions and Stirling numbers and applications}, J. Comput. Appl. Math. \textbf{260} (2014), 201\nobreakdash--207; Available online at \url{http://dx.doi.org/10.1016/j.cam.2013.09.077}.

\bibitem{exp-reciprocal-cm-IJOPCM.tex}
X.-J. Zhang, F. Qi, and W.-H. Li, \textit{Properties of three functions relating to the exponential function and the existence of partitions of unity}, Int. J. Open Probl. Comput. Sci. Math. \textbf{5} (2012), no.~3, 122\nobreakdash--127.

\end{thebibliography}
\end{document}